\documentclass{pnastwo}
\usepackage{mathrsfs}
\usepackage{amsfonts}

\usepackage[dvips]{graphicx}
\usepackage{pnastwoF}
\usepackage{amsmath,amssymb,amsfonts,eepic,epic,color}

\url{www.pnas.org/cgi/doi/10.1073/pnas.0709640104}
\copyrightyear{2008}
\issuedate{Issue Date}
\volume{Volume}
\issuenumber{Issue Number}

\begin{document}

\title{The evolution of the random reversal graph}

\author{Christian M. Reidys\affil{1}
       {Center for Combinatorics, LPMC-TJKLC\\Nankai University \\
         Tianjin 300071\\
         P.R.~China\\
         Phone: *86-22-2350-6800\\
         Fax:   *86-22-2350-9272\\
         duck@santafe.edu}
\and Emma Y. Jin}

\contributor{Submitted to Proceedings of the National Academy of
Sciences of the United States of America}

\maketitle

\begin{article}
\begin{abstract}
Although genomes can in terms of letter-by-letter DNA nucleotide
content be very similar, it is possible that they differ
significantly by gene order as well as gene orientation. Order and
orientation of genes is represented as signed permutations and we
study here genome rearrangements via reversals (Sankoff et al.
(1999) Proceedings of the third annual international conference on
Computational molecular biology, 302-309), flipping entire segments
of genes and changing their orientations (signs). This abstraction
leads to the notion of the reversal graph in which two signed
permutations are neighbors if they differ by one reversal. The
structure of this graph is consequently of importance for
understanding the evolution of genomes. In this paper we identify
the random reversal graph, i.e.\ the probability space consisting of
subgraphs over signed permutations, obtained by selecting edges
(reversals) with independent probability $\lambda_n$ as an important
framework for the analysis of genome rearrangements. One key control
parameter of the reversal graph, the rate of reversals, has already
been quantified using comparative analysis (Seoighe, et al. (2000)
Proc. Natl. Acad. Sci. USA 97: 14433-14437) (Sharakhov, et al.
(2002) Science 298: 182-185). The random reversal graph offers new
perspectives, allowing to study the connectivity of genomes as well
as their most likely distance as a function of the reversal rate.
Our main result shows that the structure of the random reversal
graph changes dramatically at $\lambda_n=1/\binom{n+1}{2}$. For
$\lambda_n=(1-\epsilon)/\binom{n+1}{2}$, the random graph consists
of components of size at most $O(n\ln(n))$ a.s.\ and for
$(1+\epsilon)/\binom{n+1}{2}$, there emerges a unique largest
component of size $\sim \wp(\epsilon) \cdot 2^n\cdot n!$ a.s.. This
``giant'' component is furthermore dense in the reversal graph.
\end{abstract}


\abbreviations{}

\section{Introduction}

\dropcap{S}eventy years ago, Dobzhansky and Sturtevant
\cite{Dobzhansky} initiated the study of genome rearrangements. Due
to recent progress in large-scale sequencing and comparative
mapping, genome rearrangements have become increasingly popular
\cite{Brien,Murphy,Lander,Venter}. One central question in this
context is that of the true evolutionary distance between two
genomes\footnote{We remark here upon the notion of multichromosomal
rearrangement problems, where plausible rearrangement scenarios for
multiple species \cite{Max:07} are studied}. Unichromosome genomes
evolve by local mutations, that is insertions, deletions and
substitution of nucleotides as well as global rearrangements, i.e.\
inversions. The latter are a quite common mode of molecular
evolution in both, unichromosome as well as multichromosome genomes,
like mitochondrial, chloroplast, viral and bacterial DNA
\cite{Pevzner:95,Bafna:95a,Bafna:95b,Franklin:05,Bafna:96}. Data
show that the genomes of different species, most of which being very
similar on the letter-by-letter DNA code basis, differ mainly in
gene order and orientation.

For instance, the mitochondrial genomes of Brassica (cabbage) and
Brassica campestris (turnip) are closely related having more than
$99\%$ identical genes \cite{Palmer:88}. Remarkably, mouse and human
carry about $30,000$ genes, which are, on DNA letter comparison
basis, to $85\%$ identical. Nevertheless, the significant
differences between the two species manifest in terms of order and
orientation of the four-letter DNA codes. In fact, comparing the two
genomes provides an evolutionary history of the two species and
traces out the diverging sequences of DNA.

The above observations motivated Kececioglu, Sankoff and others
\cite{Sankoff:92} to model the evolution of genomes by representing
order and orientation of genes as signed permutations, see Fig.\
\ref{F:humanworm}.

For uni-chromosome genomes we have according to
\cite{Sankoff:96,Wang}, the following rearrangement mechanisms {\bf
(a)} reversals (inversions), of any number of consecutive terms in
the ordered set, which, in case of signed orders also reverses the
polarity of each term within the scope of the inversion.
\cite{Sankoff,Pevzner:95,Sankoff:96} and {\bf (b)} transpositions of
any number of consecutive terms from their position in the order to
a new position between any other pair of consecutive genes. In
addition, we have in multi-chromosome genomes {\bf (c)} reciprocal
translocations.

While it is plain that there are multiple edit operations acting on
the genome, it is still unclear how to weight them. For instance,
transpositions are observed much less frequently than reversals
\cite{Sharakhov:02}. In this paper we restrict our analysis to
reversals as edit operation and study monotone properties, that is
graph properties that remain valid when adding more
edges\footnote{for instance ``connectivity'' is a monotone property:
a graph remains connected when additional edges are inserted}.
Studying monotone properties allows us to derive conclusions that
remain valid even if additional rearrangement operations are taken
into account.

Reversals change the order of the genes in a genome, and also the direction of
transcription. For instance, the reversal $\rho_{3,5}$ maps the
signed permutation $(+1,+4,\boxed{+2,+5,+3})$ to
$(+1,+4,\boxed{-3,-5,-2})$. Reversals and signed permutations
constitute the reversal graph in which two signed permutations are
adjacent if they differ by a single reversal.
The reversal graph is regular of degree $\binom{n+1}{2}$, has
diameter $(n+1)$ \cite{diam} and is connected, since any signed
permutation can be generated via a sequence of reversals.

The reversal graph implies an obvious notion of distance, that is the
minimal number of reversals needed to transform one genome into the other.
Computing this reversal distance turns out to be a difficult problem.
For (unsigned) permutations with only reversal moves allowed, it has been
shown to be NP-hard \cite{Caprara:97}. For signed permutations,
Hannenhalli and Pevzner \cite{Pevzner:95} did breakthrough work and
derived an algorithm capable of computing it in polynomial time.
Since then, Pevzner {\it et al.}\ \cite{Pevzner:03a,Pevzner:03b},
Gog {\it et al.}\ \cite{Bader:08}, Kaplan {\it et al.}\ \cite{Kaplan:96}
and Berard {\it et al.}\ \cite{Berard:07} continued to work on the
parsimonious rearrangement scenario.

A question of central interest is that of the true evolutionary
distance (TED) between two genomes. The TED has been studied by
means of random reversal walks \cite{Caprara}, where it is
oftentimes related to the reversal distance \cite{Wang,
Caprara,Pevzner,Berestycki:06}. Also Markov chain Monte Carlo (MCMC)
models have been used in order to identify evolutionary trajectories
\cite{Miklos,Miklos:03,Larget}. A random reversal walk is a walk
that starts at the identity and is obtained by reversal moves,
selecting at each node incident reversals uniformly at random. It is
well-known that for distant genomes, the TED is generally much
greater than the shortest distance
\cite{Wang,Caprara,Pevzner,Berestycki:06}. Accordingly, the actual
set of paths connecting two signed permutations is only implicitly
considered. There is, for instance, no concept of the distribution
of the lengths of such paths.

Since the computation of random reversal distances is difficult,
models better suited for computation are shifted into focus. As a
result, random transposition walks in the symmetric group were
studied \cite{Berestycki:06}, where the transposition distance of a
random transposition walk after $t$ steps is analyzed. In
\cite{Berestycki:06} the authors find that in the supercritical
regime results of random transposition walks also apply to random
reversal walks, since fragmentation of the breakpoint graph induced
by reversals can be ignored \cite{Berestycki:06}. Eriksen {\it et
al.} \cite{Eriksen} also study random transposition walks and derive
a closed formula for the expected distance by considering Markov
chains on the conjugacy classes of the permutation group. It is then
shown in \cite{Eriksen} how to approximate the expected reversal
distance.

In the following, we analyze the structure of the reversal graph in
the language of random graph theory. This theory has been founded by
Paul Erd\" os and Alfred Ren\'yi in 1960 \cite{Erdos,Erdos:61}.
Erd\" os and Ren\'yi discovered a rather dramatic structural change
of random subgraphs of the complete graph obtained by selecting
edges with probability $\lambda_n$, around $\lambda_n=1/n$. Let
$\lambda_n=(1-\epsilon)/n$ where $\epsilon$ is some small, positive
constant. Then the random graph consists of small components, the
largest of which having size $O(\ln(n))$. In case of
$\lambda_n=(1+\epsilon)/n$, however, the scenario changes: many of
the small components have merged, forming a unique, giant component
of size $O(n)$. The remaining vertices remain organized in small
components, having size at most $O(\ln(n))$.

Here we follow the random graph paradigm and study random subgraphs
of the reversal graph, obtained by selecting reversals with
probability $\lambda_n$, having all signed permutations as vertices, see
Fig.~\ref{F:para-c1}.
In analogy to Erd\" os and Ren\'yi's result, we also observe a phase
transition at $\lambda_n^*=1/\binom{n+1}{2}$. Selecting reversals
with probability $(1-\epsilon)/\binom{n+1}{2}$ the random graph
consists of small components. Increasing this probability to
$(1+\epsilon)/\binom{n+1}{2}$, these components merge into a unique
giant of size precisely $\sim\wp(\epsilon)\cdot 2^n\cdot n!$, where
$\wp(\epsilon)$ denotes the survival probability of a certain
branching process.

The paper is organized as follows: we first provide some necessary
background. Then we present the main result and discuss the key
ideas of its proof. Detailed arguments of all technical lemmas are
presented in the Supplemental Materials (SM).


\section{Random graphs}


\subsection{Some basic facts}
A graph $G$ is a pair consisting of a vertex set and an edge set.
We will identify edges with pairs of vertices and two vertices
``connected'' by an edge are called adjacent, or neighbors.
We denote the size of the vertex set by $\vert G\vert$.
Suppose $G$ is a finite group and $S\subset G$ such that
$\langle S \rangle=G$, $S=S^{-1}$ and $1\not \in S$. Then the Cayley
graph $\Gamma(G,S)$, has $G$ as vertex set and edges of the form
$\{v,v'\}$, where $v^{-1}v'\in S$.

For $v,v'\in G$, let $d(v,v')$ be the length of the shortest
$G$-path connecting $v$ and $v'$ or $\infty$, otherwise. For
$A\subset G$, we set $\text{\sf B}(A,j)$ to be the set of all
elements of distance smaller or equal to $j$ to elements of $A$
and $\text{\sf d}(A)$ to be the set of elements, not contained in $A$,
that are adjacent to some $a\in A$. We call $\text{\sf B}(A,j)$ and
$\text{\sf d}(A)$ the ball of radius $j$ around $A$ and the vertex
boundary of $A$ in $G$, respectively.
If $A$ contains only one element, $\alpha$, we simply write $\text{\sf
B}(\alpha,j)$. A subset $E\subset G$, is called dense in $G$ if
$\text{\sf B}(g,1)\cap E\neq\varnothing$ for any $g\in G$.

Given a graph $G_n$, the random graph is the probability space
consisting of $G_n$-subgraphs, obtained by independently selecting
edges\footnote{Plainly we can also consider induced subgraphs via
selecting vertices} with probability $\lambda_n$. Since all sets
involved are finite, taking the power set of subgraphs of $G_n$ as
Borel $\sigma$-algebra, the random graph $W$ is tantamount to the
set of subgraphs and the probability
\begin{equation}\label{E:proba}
\mathbb{P}_n(W)=\lambda_n^{e_W}(1-\lambda_n)^{e_{G_n}-e_W},
\end{equation}
where $e_{G_n}$ and $e_W$ denote the number of edges contained in $G_n$ and
$W$, respectively.

A property, $\mathcal{P}$, is tantamount to a set of subgraphs
closed under graph isomorphisms. $\mathcal{P}$ is called monotone,
if increasing $\lambda_n$ increases $\mathbb{P}_n(\mathcal{P})$
i.e.\ if the mapping $\lambda_n\mapsto\mathbb{P}_n(\mathcal{P})$ is
monotone. It is known \cite{Bollobas:86} that monotone properties
exhibit thresholds. That is, there exists a critical probability
$\lambda_n^*$ such that, for large $n$, we have the following
situation:
\begin{equation*}
\mathbb{P}_n(\mathcal{P})=
\begin{cases}
o(1)      & \text{\rm for $\lambda_n=o(\lambda_n^*)$}\\
1-o(1)    & \text{\rm for $\lambda_n^*=o(\lambda_n)$.}
\end{cases}
\end{equation*}

A component is a maximal, connected subgraph, $C_n$, and the
largest component of the random graph is denoted by $C_n^{1}$.
It is called a giant if and only if any other component satisfies
$\vert C_n\vert = o(\vert C^1_n\vert)$.

\subsection{The random reversal graph}

Let $S_n$ denote the symmetric group over $[n]$. A permutation
is a one-to-one mapping from $x\colon [n]\longrightarrow [n]$
and represented as an $n$-tuple
$
x=(x_1,\dots,x_n)
$,
where $x_i=x(i)$. Furthermore, let
$\varepsilon=(\varepsilon_1,\dots,\varepsilon_n) \in \{-1,+1\}^n$
denote the $n$-tuple of ``signs''.

A signed permutation is a pair $(\varepsilon,x)=(\varepsilon_1 x_1,
\dots,\varepsilon_nx_n)$ and we denote the set of signed
permutations by $B_n$. Plainly we have $\vert B_n\vert =2^nn!$.
$B_n$ carries a natural structure of a group via
\begin{equation*}
(\varepsilon_x,x)\cdot (\varepsilon_y,y)=
(\varepsilon_x\cdot \varepsilon_y^{x},x\cdot y)\quad
\text{\rm where}\ (\varepsilon_y^{x})_i=(\varepsilon_y)_{x(i)}.
\end{equation*}
That is, there is an additional action on $\varepsilon_y$ when
commuting it with $x$, given by the $x$-permutation of the coordinates.
A reversal $\rho_{i,j}$ is a special type of signed permutation
\begin{equation*}
\rho_{i,j} = (\varepsilon^\rho_{i,j},(1,\dots, i-1, j, j-1, \dots,
i, j+1,\dots))
\end{equation*}
where
\begin{equation*}
(\varepsilon^\rho_{i,j})_h=
\begin{cases}
-1 & \text{\rm for $i\le h\le j$}\\
+1 & \text{\rm otherwise.}
\end{cases}
\end{equation*}
Accordingly, a reversal $\rho_{i,j}$ acts via right-multiplication
in $B_n$ as follows
\begin{eqnarray*}
(\varepsilon_1x_1,\dots,\varepsilon_ix_i,\dots,\varepsilon_jx_j,
\dots,\varepsilon_nx_n)\cdot \rho_{i,j} = \\
(\varepsilon_1x_1,\dots,-\varepsilon_jx_j,\dots,-\varepsilon_ix_i,
\dots, \varepsilon_nx_n).
\end{eqnarray*}
In other words $\rho_{i,j}$ transforms the subsequence
$(\varepsilon_ix_i,\dots,\varepsilon_j x_j)$ into
$(-\varepsilon_jx_j,\dots,-\varepsilon_ix_i)$ by inverting the order
and the signs within the interval $[i,j]$.
For instance
\begin{eqnarray*}
(+5,+2,-1,+3,-4)\cdot \rho_{2,3} & = &
(+5,+1,-2,+3,-4) \\
(+5,+2,-1,+3,-4)\cdot \rho_{2,2} & = & (+5,-2,-1,+3,-4)
\end{eqnarray*}
and we notice that one can change the sign of individual coordinates
via the reversals $\rho_{i,i}$. Let $R_n\subset B_n$
denote the set of all reversals and let $\Gamma(B_n,R_n)$ denote the
reversal graph.  We observe that the reversal graph is
regular of degree $\binom{n-1}{2}+n=\binom{n+1}{2}$ and has diameter
$(n+1)$ \cite{diam}. It is furthermore connected, since any
signed permutation can be generated starting with the identity via a
sequence of reversals. A $s$-cell is a connected subgraph of the reversal graph
of size at least $O(n^s)$ and we set $\Delta_n=\Gamma(B_n,R_n)$.

Having constructed the reversal graph, we are now in position to introduce
the random reversal graph as a probability space consisting of
subgraphs of the reversal graph. These subgraphs, $\Gamma_n$, have vertex
set $B_n$ and are obtained by selecting each edge with probability
$\lambda_n$. Following the logic of eq.\ {\bf [\ref{E:proba}]} a specific
subgraph has the probability
\begin{equation*}
\mathbb{P}_n(\Gamma_n)=\lambda_n^{e_{\Gamma_n}}(1-\lambda_n)^{
\frac{1}{2}\binom{n+1}{2}2^nn!-e_{\Gamma_n}},
\end{equation*}
where $e_{\Gamma_n}$ denotes the number of edges contained in $\Gamma_n$.


\section{The main result}

Let $\epsilon$ denote some positive constant smaller than one and
$n^{-\frac{1}{4}+\delta}\le \epsilon_n<1$ for some $0<\delta<\frac{1}{4}$.
Suppose further $0<x<1$ is the unique root of $e^{-(1+\epsilon)y}=1-y$, and
\begin{equation}\label{E:it}
\wp(\epsilon_n)=
\begin{cases}(1+o(1)) x & \text{\rm for $\epsilon_n=\epsilon>0$} \\
(2+o(1))\epsilon_n & \text{\rm for
$n^{-\frac{1}{4}+\delta}\le \epsilon_n=o(1)$.}
\end{cases}
\end{equation}
We are now in position to state our main result
\par\medskip
\begin{theorem}\label{T:main}
Let $\Gamma_n$ be a random subgraph of $\Delta_n$,
obtained by selecting edges with probability $\lambda_n$. Then,
a.s.\\
\begin{equation*}
\vert C_n^{1}\vert \sim
\begin{cases}
O(n\ln(n)) & \quad \text{\it for }
\lambda_n=(1-\epsilon)/\binom{n+1}{2} \\
\wp(\epsilon_n)\ 2^n\cdot n! & \quad \text{\it for }
\lambda_n=(1+\epsilon_n)/\binom{n+1}{2}
\end{cases}
\end{equation*}
and for $\lambda_n=(1+\epsilon_n)/\binom{n+1}{2}$, the largest component
$C_n^{1}$ is a.s.\ unique.
\end{theorem}

In order to prove Theorem~\ref{T:main}, we establish several auxiliary results.
The first two of these are concerned with the formation of $k\delta$-cells.
Note that since $\delta>0$ is fixed, choosing $k$ sufficiently large,
these cells have arbitrarily high polynomial degree.

The first step for generating $k\delta$-cells is the simulation of a
branching process in $\Delta_n$. The purpose of this process is
the generation of an acyclic, connected subgraph of size $\lfloor
\frac{1}{4}n^{\frac{3}{4}}\rfloor$.

\begin{lemma}\label{L:cell1}
Each signed permutation $v$ is contained in a connected, acyclic
subgraph, $\mathcal{T}_n(v)$ of size
$\lfloor\frac{1}{4}n^{\frac{3}{4}} \rfloor$ with probability at
least $\wp(\epsilon_n)$ given via eq.~{\bf [\ref{E:it}]}.
\end{lemma}

The second step consists in showing that the graphs $\mathcal{T}_n(v)$
generated in Lemma~\ref{L:cell1} are the building blocks of
$k\delta$-cells.
This works inductively using sets of reversals over coordinates
that were not considered during the generation of $\mathcal{T}_n(v)$.
The intuition here is to view the coordinates of a permutation as
``dimensions'' in which specific reversals act. Careful partitioning
of these coordinates eventually allows to merge $\mathcal{T}_n(v)$-subgraphs
to $k\delta$-cells, see Fig.\ \ref{F:cell2}.
\begin{lemma}\label{L:cell2}
Suppose $k\in\mathbb{N}$ is arbitrary but fixed. Then
each $\Gamma_n$-vertex is contained in a $k\delta$-cell with
probability at least
\begin{equation*}
\delta_k(\epsilon_n)=\wp(\epsilon_n)\, (1-e^{-\beta_k\theta_n}),
\quad \text{\it for some }\beta_k>0.
\end{equation*}
\end{lemma}

Let $\Gamma_{n,k}$ denote the set of vertices contained in
$k\delta$-cells. We next show that the number of these vertices
is sharply concentrated.

\begin{lemma}\label{L:small1}
Suppose $k\in\mathbb{N}$ is sufficiently large. Then
\begin{equation*}
\vert \Gamma_{n,k}\vert \sim \wp(\epsilon_n)\cdot 2^n\cdot n! \
\qquad \text{\it a.s.~.}
\end{equation*}
\end{lemma}

Before we can prove the main theorem we observe that
$\Gamma_{n,k}$ is dense in $B_n$ a.s. (SM, Lemma~$7$). Using this
density we show in Lemma~\ref{L:split} below that there exist many, short
vertex disjoint paths between certain $\Gamma_{n,k}$-splits.

\begin{lemma}\label{L:split}
Let $(S,T)$ be a vertex-split of $\Gamma_{n,k}$ such that
there exist $0< \rho_0\le \rho_1<1$ and
\begin{equation*}\label{E:prop}
\begin{split}
2^n\cdot(n-2)!\le
\rho_0 \vert \Gamma_{n,k}\vert = \vert S\vert \le \vert T \vert =
\rho_1 \vert \Gamma_{n,k}\vert .
\end{split}
\end{equation*}
Then there exists some $c>0$ such that $S$ is connected to $T$ in
$\Delta_n$ via at least
\begin{equation*}
c\cdot\frac{2^n\cdot(n-3)!}{\binom{n+1}{2}^3}
\end{equation*}
edge disjoint (independent) paths of length $\le 3$, a.s..
\end{lemma}
We remark that Lemma~\ref{L:split} does not
use an isoperimetric inequality \cite{Harper:66b}. It only employs a
generic estimate of vertex boundaries in Cayley graphs due to Aldous
\cite{Aldous:87} and Babai \cite{Babai:91b}.

The idea here is the following: if $\Gamma_{n,k}$ does not contain
the giant component, then there exists a certain type of split. The
paths of Lemma~\ref{L:split} would then connect the two subsets of
this split and none of them can be chosen
\cite{Ajtai:82,Reidys:08rand} in the random graph. To prove the main
result, our strategy will be to show that the probability of
choosing none of these paths tends to zero.
\section{Proof of Theorem~\ref{T:main}.}

To prove the theorem, we follow Ajtai {\it et al.}
\cite{Ajtai:82,Reidys:08rand} and select the $\Delta_n$-edges in two
distinct randomizations. Suppose $x_1,x_2>0$ such that
$x_1^{-1}+x_2^{-1}=1$. First we select with probability
$(1+\epsilon_n/x_1)/\binom{n+1}{2}$ and second with probability
$(\epsilon_n/x_2)/\binom{n+1}{2}$. The probability of not being
chosen in both rounds is given by
\begin{equation*}
\left(1-\frac{1+\epsilon_n/x_1}{\binom{n+1}{2}}\right)
\left(1-\frac{\epsilon_n}{x_2\cdot \binom{n+1}{2}}\right)\ge
1-\frac{1+\epsilon_n}{\binom{n+1}{2}},
\end{equation*}
whence it suffices to prove that after the second randomization,
there exists a giant component with the property $\vert
C_{n}^{1}\vert\sim \vert\Gamma_{n,k}\vert$. After the first
randomization each $\Delta_n$-edge has been selected with
probability $(1+\epsilon_n/x_1)/\binom{n+1}{2}$ and according to
Lemma~\ref{L:small1}, we have
\begin{equation*}
\vert \Gamma_{n,k}(x_1)\vert\sim \wp(\epsilon_n/x_1)\cdot 2^n\cdot
n!\quad \text{\rm a.s.}
\end{equation*}
Suppose $\Gamma_{n,k}(x_1)$ contains a split, that is
two ``large'' components, $S,T$. Assume $\vert S\vert\le
\vert T\vert$ where
$$
(n-2)!\cdot 2^n \le \vert S\vert \le (1-b)\, \vert
\Gamma_{n,k}(x_1)\vert, \quad \text{\rm where } b>0.
$$
Then there exists a split of $\Gamma_{n,k}(x_1)$, denoted by
$(S,T)$, satisfying the assumptions of Lemma~\ref{L:split} (and, of
course, $S\cap T=\varnothing$). We observe that Lemma~\ref{L:cell2}
limits the number of ways these splits can be constructed. Suppose
$M_k\ge c_k \, n^{k\delta+\frac{3}{4}}$ for some $c_k>0$. Obviously,
there are at most $2^{2^n\cdot n!/M_k}$ ways to select $S$ in such a
split. Now we employ Lemma~\ref{L:split}. In view of
$2^n\cdot(n-2)!\le \vert S\vert$, Lemma~\ref{L:split} implies that
there exists some $c>0$ such that $S$ is connected to $T$ in
$\Delta_n$ via at least $c\cdot 2^n(n-3)!/\binom{n+1}{2}^{3}$ edge
disjoint paths of length $\le 3$, a.s.. We next perform the second
randomization and select $\Delta_n$-edges with probability
$(\epsilon_n/x_2)/\binom{n+1}{2}$. None of the above paths can be
selected during this process. Since any two paths are edge disjoint,
the expected number of such splits is, by linearity of expectation,
less than
\begin{equation*}
\begin{split}
2^{2^n\cdot n!/M_k}
\left[1-\left(\frac{\epsilon_n/x_2}{\binom{n+1}{2}}
\right)^3\right]^{\frac{c\cdot 2^n
(n-3)!}{n^{6}}} \le \qquad \\
2^{2^nn!/M_k} e^{-c' 2^nn!/n^{15}}
\quad\text{for some $c,c'>0$}.
\end{split}
\end{equation*}
Accordingly, choosing $k$ sufficiently large, the expected number of
these $\Gamma_{n,k}(x_1)$-splits tends to zero, i.e.~for any $k\ge
k_0\in \mathbb{N}$, there exists no two component split $(S,T)$ of
$\Gamma_{n,k}(x_1)$ with the property
$\rho_0\vert\Gamma_{n,k}(x_1)\vert=\vert S\vert\le \vert T\vert$
a.s.. Consequently, there exists some subcomponent $C_n(x_1)$ with
the property
\begin{equation*}
\vert C_n(x_1) \vert=\vert \Gamma_{n,k}(x_1)\vert \sim
\wp(\epsilon_n/x_1)\, \cdot 2^n\cdot n! \quad \text{\rm a.s.,}
\end{equation*}
obtained by the merging of the $k$-cells generated during the first
randomization via the paths selected during the second. Since
$\wp(\epsilon_n/x_1)$ is continuous in the parameter $\epsilon_n/x_1$,
see eq.\ {\bf [\ref{E:it}]}, we derive, for $x_1$ tending to $1$
\begin{equation*}
\vert C_n^{1} \vert=\lim_{x_1\to 1}\vert C_n(x_1) \vert \sim
\wp(\epsilon_n) \cdot 2^n\cdot n!\quad \text{\rm a.s.}
\end{equation*}
It remains to prove uniqueness. Obviously, any other largest
component, $\tilde{C}_n$, is necessarily contained in
$\Gamma_{n,k}$. However, we have just proved $\vert C_n^{1} \vert
\sim \wp(\epsilon_n)\cdot 2^n\cdot n!$ and according to
Lemma~\ref{L:small1}, $\wp(\epsilon_n)\cdot 2^n\cdot n!\sim \vert
\Gamma_{n,k}\vert$. Therefore $\vert\tilde{C}_n\vert=o(\vert
C_n^{1}\vert)$, whence $C_n^{1}$ is unique. \hfill $\square$



\section{Discussion}


In \cite{Pevzner,Berestycki:06} we find general consensus that, if the number
of reversals is smaller than $0.4n$, then the parsimonious distance between
genomes of length $n$ approximates the true evolutionary distance well.
In this context Pevzner proposed the hypothesis of {\it fragile sites},
i.e.\ sites, that are acted upon by reversals with higher probability.
While this hypothesis implies that the parsimonious distance
is a key notion for true evolutionary distances, it does not shed light
on the mechanisms of genome rearrangements.

Even for metazoa mitochondrial genomes (MMGs), that contain only
single copies of genes and whose rearrangements lack transposons as
well as introns that potentially affect mutations \cite{Miklos}
reversals are not the only rearrangement mechanism
\cite{Sankoff:96,Wang,Bernt,Perseke}. Nevertheless they represent an
important edit operation and since MMGs exhibit no preference of a
particular ordering of the genes \cite{Miklos}, they are
particularly well suited to be modeled in terms of signed
permutations and reversals.

Random reversals are not a new concept. They constitute for instance the basic moves
of random reversal walks \cite{Caprara}. Although there are obvious connections
between these walks and the random reversal graph both approaches complement
each other. For instance, to our knowledge, within the framework
of random reversal walks, there is no notion of random variables as a function
of the rate of reversals.

In \cite{Fischer:06} reversal rates are quantified with respect to synteny
blocks (s-blocks), i.e.\ chromosomal regions in different genomes sharing a
common evolutionary origin. In other words, two chromosomal regions are
syntenic when multiple consecutive genes are found in a conserved
order between the two genomes \cite{Fischer:06}. In Table\ 1 \cite{Fischer:06}
we present reversal rates within
synteny blocks between the species of Hemiascomycetes. The table displays
the mean number of genes per block, the mean number of inversions
per block and the mean frequency of gene inversions between the
respective pairs of species from different lineages.
{\small
\begin{table}
\begin{tabular}{l|c|c|c|c}
Species & genes/s-block  & rev/s-block & rev-rate\\
\hline
 S.cerevisiae-C.glabrata     &  7.6   & 0.2 & 0.08\\
 S.cerevisiae-K.lactis       &  6.0   &  0.2  & 0.05\\
 S.cerevisiae-A.gossypii     &  6.6   &  0.3  & 0.09 \\
 C.glabrata-K.lactis      &  5.0   &  0.2  & 0.07 \\
 C.glabrata-A.gossypii     &  5.0   &  0.2  & 0.08 \\
 K.lactis-A.gossypii      &  8.8   &  0.1  & 0.03\\
\hline
 S.cerevisiae-D.hansenii     &  2.6   & 1.1  & 0.57 \\
 C.glabrata-D.hansenii       &  2.5   &  1   & 0.52 \\
 K.lactis-D.hansenii         &  2.7   &  1.3 & 0.65 \\
 A.gossypii-D.hansenii       &  2.9   &  1.3 & 0.60  \\
\hline
\end{tabular}
{\small {\bf Tab.\ 1}:comparative analysis of reversals within synteny
blocks of hemiascomycete yeasts, adapted from \cite{Fischer:06}.}
\end{table}
}
According to Table\ 1, the mean number of genes per block ranges
from $2.5$ to $8.8$, implying that the inversions are of short and the
frequency of gene inversion ranges from $0.03$ to $0.65$. In the random reversal
model the critical probability for the emergence of the giant component is
given by $1/\binom{n+1}{2}$. Accordingly, for the genomes of Table\ 1,
ranging in length from $2.5$ to $8.8$, the
corresponding critical probability ranges from $0.02$ to $0.23$.
Thus the reversal rates of Table\ 1 are above the
threshold value for the existence of a giant component.
We finally remark, that Table\ 1 is derived from the algorithm
{\tt DERANGE} \cite{Sankoff:96}, which is based on the parsimonious
approach.

Further studies on the rate of reversals were obtained by Sharakhov
{\it et al.} \cite{Sharakhov:02}, who studied gene reversals in
Anopheles gambiae and A.\ funestus. The authors estimated the
reversal rate on the basis of time, considering the mean length of
conserved segments in which the gene order are preserved during
mutation \cite{Nadeau}. Also Seoighe {\it et al.} \cite{Seoighe}
estimate the rates of small reversals relative to that of large
rearrangements, including translocations, large reversals and long
distance transpositions by means of comparative analysis of the
species C.\ albicans and S.\ cerevisiae.

Plainly, if paths between two genomes are a result of random reversals,
they are not necessarily of minimal length. Furthermore, the length
of these reversal paths may not be the most important question. For
instance, it could be much more relevant to identify the most ``likely''
path between two genomes.

The random reversal graph introduced here allows us to study the
above question and augments the probabilistic analysis of random
reversal walks. The fragile site hypothesis has a natural
interpretation within our framework and we can incorporate the
probability of reversals as model parameter. As for the former, a
weighted random graph can be studied, in which reversals are no
longer uniformly chosen. Similar ideas exist already in random graph
theory, for instance in the context of random graphs with given
average degree sequence \cite{Fan}. As for the latter, our main
result shows this as follows: before the phase transition,
the random reversal graph contains only very small components. Thus
a reversal path between two genomes is essentially non-existent,
leading to a path length of infinity. However, after the phase
transition, two genomes are connected via a reversal path with some
positive probability, i.e.\ the length of a reversal path becomes
finite.

While it is beyond the scope of this paper to provide the detailed
analysis of the length of these paths, our framework enables us to
study the following questions:\\
{\bf (a) } Given two signed permutations, $\alpha$, $\beta$ compute the
threshold probability, $\lambda_n^\sharp$ for which there exists
some $c_{\alpha,\beta}>0$ such that
\begin{equation}\label{E:mono}
d_{\Gamma_n}(\alpha,\beta) \le c_{\alpha,\beta}\cdot
                                  d_{\Delta_n} (\alpha,\beta)
\quad \text{\rm a.s.,}
\end{equation}
where $d_{\Gamma_n}$ and $d_{\Delta_n}$ denote the distance in
$\Gamma_n$ and $\Delta_n$, respectively. Indeed $\lambda_n^\sharp$
is guaranteed to exist by general theory \cite{Bollobas:86}, since
eq.\ {\bf [\ref{E:mono}]} represents a
monotone property of the random reversal graph.\\
{\bf (b)} how many independent
$\Gamma_n$-paths of length at most
$c_{\alpha,\beta}\cdot d_{\Delta_n}(\alpha,\beta)$ exist?\\
As for {\bf (a)} the critical probability is
$\lambda^{\sharp}=O\left(1/n\right)$
and that these paths can be constructively generated following the
ideas in \cite{Reidys}. The intuition here is that given two signed
permutations, random reversal paths between them can be derived via
a ``parallel'' branching, essentially by shifting minimal reversal
paths between them. One can
paraphrase the situation by saying that if individual genes have for
sufficiently long genomes, a positive probability of being
rearranged, then the minimal reversal distance approximates up to a
constant factor the true evolutionary distance. Furthermore, this
constant factor can be explicitly calculated. Assertion {\bf (b)}
provides a type of information previously unavailable regarding the
variety of reversal paths.

We remark that the ideas in this paper show why the simulation of
reversal walks via sign-change transposition models
\cite{Berestycki:06,Eriksen} works so well. These models assume
``simpler'', transposition based rearrangement moves and have proven
to approximate the more complex reversal walks with high accuracy
\cite{Berestycki:06}.

Indeed, probabilistic analysis reveals that the random transposition
graph and the random reversal graph are structurally very similar.
In order to elaborate on this further, we introduce a sign-change
transposition as a particular type of signed permutation
\begin{equation*}
\tau_{i,j}= (\varepsilon^\tau_{i,j},
(1,\dots,i-1,j,i+1,\dots,j-1,i,j+1,\dots));\
\end{equation*}
where
\begin{equation*}
(\varepsilon^\tau_{i,j})_h=
\begin{cases}
-1 & \text{\rm for $h=i,j$}\\
+1 & \text{\rm otherwise.}
\end{cases}
\end{equation*}
Accordingly, a sign-change transposition, $\tau_{i,j}$ acts by
mapping $\varepsilon_ix_i$ into $-\varepsilon_jx_j$ and
$\varepsilon_jx_j$ into $-\varepsilon_ix_i$ and as identity,
otherwise. Following the ideas in this paper it can be shown
that the following analogue of Theorem~\ref{T:main} holds:
\begin{equation*}
\vert C_n^{1}\vert \sim
\begin{cases}
O(n\ln(n)) & \quad \text{\it for }
\lambda_n=(1-\epsilon)/\binom{n+1}{2} \\
\wp(\epsilon_n)\ 2^n\cdot n! & \quad \text{\it for }
\lambda_n=(1+\epsilon)/\binom{n+1}{2}
\end{cases}
\end{equation*}
and for $\lambda_n=(1+\epsilon)/\binom{n+1}{2}$, the largest
component $C_n^{1}$ is a.s.\ unique.


\begin{acknowledgments}
We want to thank Markus Nebel, Peter Stadler for comments. This work
was supported by the 973 Project, the PCSIRT Project of the Ministry
of Education, and the National Science Foundation of China.
\end{acknowledgments}

\end{article}


\end{document}


\title[Supplemental Material]
      {Supplemental Material: \\
The evolution of the random reversal graph}
\author{Christian M. Reidys$^{\,\star}$ and Emma Y. Jin}
\address{Center for Combinatorics, LPMC-TJKLC \\
         Nankai University  \\
         Tianjin 300071\\
         P.R.~China\\
         Phone: *86-22-2350-6800\\
         Fax:   *86-22-2350-9272}
\email{duck@santafe.edu}
\thanks{}
\keywords{} \subjclass[2000]{05A16}
\date{November, 2009}
\maketitle {{\small
}}

\section{Some basic facts}\label{S:2}

We consider the Cayley graph $\Gamma(B_n,R_n)$, having vertex set
$B_n$ and edges $\{v,v'\}$ where
$v^{-1}v'\in R_n$. Let $B_n$ denote the set of signed permutation
of length $n$ and $R_n$ be the set of reversals
$\rho_{i,j}$ where $1\le i\le j\le n$. For $v,v'\in B_n$, let
$d(v,v')$ be the minimal number of reversals by which $v$ and $v'$
differ. For $A\subset B_n$, we set $\text{\sf B}(A,j)=\{v\in B_n\mid
\exists \, \alpha\in A;\,d(v,\alpha)\le j\}$ and $\text{\sf d}(A) =
\{v\in B_n\setminus A\mid \exists \, \alpha\in A;\, d(v,\alpha)=1\}$
and call $\text{\sf B}(A,j)$ and $\text{\sf d}(A)$ the ball of
radius $j$ around $A$ and the vertex boundary of $A$ in
$\Gamma(B_n,R_n)$. If $A=\{\alpha\}$ we simply write $\text{\sf
B}(\alpha,j)$. Let $E\subset B_n$, we call $E$ dense in $B_n$ if
$\text{\sf B}(\sigma,1)\cap E\neq\varnothing$ for any $\sigma\in
B_n$. Let ``$<$'' be the following linear order over $\Gamma(B_n,R_n)$,
$\sigma < \tau$ if and only if
$\sigma<_{\text{\rm lex}} \tau$, where $<_{\text{\rm lex}}$ denotes
the lexicographical order. Any notion of minimal or smallest element
in a subset $A\subseteq B_n$ refers
to the above linear order.

The random graph $\Gamma_{\lambda_n}(B_n,R_n)$ is the probability
space consisting of $\Gamma(B_n,R_n)$-subgraphs, $\Gamma_n$, having
vertex set $B_n$, obtained by selecting each $\Gamma(B_n,R_n)$-edge
with independent
probability $\lambda_n$. A property $\text{\sf M}$ is a subset of
induced subgraphs of $\Gamma(B_n,R_n)$ closed under graph
isomorphisms. The terminology ``$\text{\sf M}$ holds a.s.'' is
equivalent to $\lim_{n\to\infty}\text{\sf Prob}(\text{\sf M})=1$. A
component of $\Gamma_n$ is a maximal, connected, induced
$\Gamma_n$-subgraph, $C_n$. The largest $\Gamma_n$-component is
denoted by $C_n^{1}$. We write $x_n\sim y_n$ if and only if (a)
$\lim_{n\to\infty}x_n/y_n$ exists and (b)
$\lim_{n\to\infty}x_n/y_n=1$. We furthermore write $g(n)=O(f(n))$
and $g(n)=o(f(n))$ for $g(n)/f(n)\to \kappa$ as $n\to \infty$ and
$g(n)/f(n)\to 0$ as $n\to \infty$, respectively. A largest component
is called giant if it is unique in ``size'', i.e.~any other
component, $C_n$, satisfies $\vert C_n \vert=o(\vert C_n^{1}
\vert)$.

Let $Z_n=\sum_{i=1}^n \xi_i$ be a sum of mutually independent
indicator random variables (r.v.), $\xi_i$ having values in
$\{0,1\}$. Then we have Chernoff's
large deviation inequality \cite{Chernoff:52},
that is, for $\eta>0$ and
$c_\eta=\min\{-\ln(e^{\eta}[1+\eta]^{-[1+\eta]}),
\frac{\eta^2}{2}\}$
\begin{equation}\label{E:cher}
\text{\sf Prob}(\,\vert\,Z_n-\mathbb{E}[Z_n]\,\vert\,>
\eta\,\mathbb{E}[Z_n]\,) \le
        2 e^{-c_\eta \mathbb{E}[Z_n]}\, .
\end{equation}
$n$ is always assumed to be sufficiently large and $\epsilon$ is a
positive constant satisfying $0<\epsilon < 1$. We write the binomial
distribution as
$B_m(\ell,\lambda_n)=\binom{m}{\ell}\lambda_n^\ell\,
(1-\lambda_n)^{m-\ell}$.

Let us next recall some basic facts about branching processes,
$\mathcal{P}_m=\mathcal{P}_m(p)$ \cite{Harris:63,Kolchin:86}.
Suppose $\mathcal{P}_m$ is initiated at $\xi$. Let $(\xi_i^{(t)})$,
$i,t\in\mathbb{N}$ count the number of ``offspring'' of the
$i$th-``individual'' of
$(t-1)$th ``generation'', where the
r.v.~$\xi$ and $\xi_i^{(t)}$ are $B_m(\ell,p)$-distributed. Let
$\mathcal{P}_0=\mathcal{P}_0(p)$ denote the branching process for
which $\xi$ is $B_m(\ell,p)$- and all $\xi_i^{(t)}$ are
$B_{m-1}(\ell,p)$-distributed. Furthermore, let
$\mathcal{P}_P(\lambda)$, $(\lambda>0)$ denote the branching process
in which the individuals generate offsprings according to the
Poisson distribution, i.e., $\mathbb{P}(\xi_{i}^{(t)}=j)=
\frac{\lambda^j}{j!}e^{-\lambda}$. We consider the family of
r.v.~$(Z_i)_{i\in \mathbb{N}_0}$: $Z_0=1$ and $Z_{t} =
\sum_{i=1}^{Z_{t-1}}\xi_i^{(t)}$ for $t\ge 1$ and interpret $Z_t$ as
the number of individuals ``alive'' in generation $t$. Of particular
interest for us will be the limit
$\lim_{t\to\infty}\mathbb{P}(Z_t>0)$, i.e.~the probability of
infinite survival. We write
$$
\pi_0(p)=\lim_{t\to\infty}\mathbb{P}_0(Z_t>0), \
\pi_m(p)=\lim_{t\to\infty}\mathbb{P}_m(Z_t>0) \ \text{\rm and} \
\pi_P(\lambda)=\lim_{t\to\infty}\mathbb{P}_{P}(Z_t>0)
$$
for the survival probability of $\mathcal{P}_0(p)$,
$\mathcal{P}_m(p)$ and $\mathcal{P}_P(\lambda)$, respectively.

\begin{lemma}\label{L:galton}\cite{Bollobas:92}\\
$(1)$ For all $0\le p \le 1$, we have
$\pi_{m-1}(p)\le \pi_0(p) \le
\pi_m(p)$.\\
$(2)$ If $\lambda > 1$ is fixed, then $\pi_P(\lambda)$ is the unique
solution of $x + e^{
-\lambda x} = 1$ in the interval $0 < x < 1$.\\
$(3)$ Let $p = \frac{1+\epsilon_n}{m}$
and $0 <\epsilon_n =o(1)$. Then
\begin{equation*}
\pi_m(p) = \frac{2m\epsilon_n}{m-1}+ O(\epsilon_n^2).
\end{equation*}
In particular, if $r =m-s$ then
\begin{equation*}
\pi_r(p) = 2\epsilon_n + O(\epsilon_n/m) + O(s/m) + O(\epsilon_n^2);
\end{equation*}
and hence if $s = o(\epsilon_n\, m)$ then $\pi_r(p) = (1 + o(1))\pi_0(p)$.\\
\end{lemma}
Let us finally give the key facts about the relations between the survival
probabilities $\pi_0(p), \pi_m(p)$ and $\pi_P(\lambda)$:

\begin{corollary}\label{C:1}\cite{Bollobas:92}
$(1)$ If $p=\lambda/m$ where $\lambda>1$, then $\pi_0(p) =
(1 + o(1))\pi_P(\lambda)$.\\
$(2)$ Let $p=\frac{1+\epsilon_n}{m}$, where $0<\epsilon_n = o(1)$.
Then, if $r=m-s$ and
$s=o(m\epsilon_n)$,
\begin{equation*}
\pi_0(p)=(1+ o(1))\pi_r(p)=(2+o(1))\epsilon_n.
\end{equation*}
\end{corollary}
\section{$k$-cells}
\label{S:3}
In the following, we shall always assume
\begin{equation*}
\lambda_n=\lambda_n(\epsilon_n)=
\frac{1+\epsilon_n}{\binom{n+1}{2}}\quad\text{where}\
n^{-\frac{1}{4}+\delta}\le \epsilon_n<1\ \text{\rm and }
0<\delta<\frac{1}{4}
\end{equation*}
Suppose $x>0$ is the unique root of $e^{-(1+\epsilon)y}=1-y$ and
\begin{equation}\label{E:it}
\wp(\epsilon_n)=
\begin{cases}(1+o(1)) x & \text{\rm for $\epsilon_n=\epsilon>0$} \\
(2+o(1))\epsilon_n & \text{\rm for
$n^{-\frac{1}{4}+\delta}\le \epsilon_n=o(1)$.}
\end{cases}
\end{equation}
For $k\in\mathbb{N}$, we furthermore set
\begin{equation*}
\label{E:part}
\mu_n  =  \lfloor \frac{1}{2k(k+1)}n^{\frac{3}{4}}\rfloor,\quad
\ell_n = \lfloor\frac{k}{2(k+1)} n^{\frac{3}{4}}\rfloor.
\end{equation*}
\begin{lemma}\label{L:minigene}
Each signed permutation, $v$, is contained in a $\Gamma_n$-subtree
$\mathcal{T}_n(v)$ of size $\lfloor\frac{1}{4}n^{\frac{3}{4}}
\rfloor$ with probability at least $\wp(\epsilon_n)$,
given by eq.~(\ref{E:it}).
\end{lemma}
\begin{proof}
We shall construct the subtree $\mathcal{T}_n(v)$ by constructing a
branching process $\mathcal{P}_m(\lambda_n)$
\cite{Harris:63} within $\Gamma(B_n,R_n)$, initiated at $id$ where
$m$ is given by eq.~(\ref{E:cell-1m}).
The offspring of this branching process is
generated by the following set of reversals
\begin{equation*}
N=\left\{\rho_{l,r}\mid
\lfloor\frac{1}{2}n^{\frac{3}{4}}\rfloor+1\le l\le r\le
n\right\}\subset R_n.
\end{equation*}
We initiate the process as follows:
\begin{eqnarray*}
M_0=L_0 &=&\{id\}\subset B_n\\
U_0&=&\varnothing\subset N\\
D_0&=&\varnothing\subseteq\{\lfloor\frac{1}{2}
n^{\frac{3}{4}}\rfloor+1,\ldots,n\}.
\end{eqnarray*}
Suppose we are given $M_j$, $U_j\subset N$, $L_j\subset B_n$ and
$D_j$, the process stops at $j+1$ either when $L_j=\varnothing$ or
$\vert M_j\vert=\lfloor\frac{1}{4}n^{\frac{3}{4}}\rfloor$.
Otherwise, we consider the smallest element $\omega_j\in L_j$ and
connect among the smallest
\begin{equation}\label{E:cell-1m}
m= \binom{n-\lfloor\frac{1}{2}n^{\frac{3}{4}} \rfloor+1}{2}-
(n-\lfloor\frac{1}{2}n^{\frac{3}{4}}\rfloor)
\lfloor\frac{1}{4}n^{\frac{3}{4}}\rfloor
\end{equation}
$\omega_j$-neighbors, $x=\omega_j\cdot\rho_{\alpha_j, \beta_j}$. We
select with independent probability $\lambda_n(\epsilon_n)$, subject
to the conditions $\rho_{\alpha_j,\beta_j}\in N\setminus U_j$ and
$\alpha_j\not\in D_j$.
Note that if $\lfloor\frac{1}{4}n^{\frac{3}{4}}\rfloor-1$ vertices
have been connected, we have
\begin{eqnarray*}
|N|-|U_j\cup\{\rho_{\alpha_j,\beta}\mid\alpha_j\in D_j\}|
&=&|N|-|\{\rho_{\alpha_j,\beta}\mid\alpha_j\in D_j\}|\\
& \ge & \binom{n-\lfloor\frac{1}{2}n^{\frac{3}{4}}\rfloor+1}{2}-
(n-\lfloor\frac{1}{2}n^{\frac{3}{4}}\rfloor)
\lfloor\frac{1}{4}n^{\frac{3}{4}}\rfloor=m
\end{eqnarray*}
Therefore, as long as we connect less than
$\lfloor\frac{1}{4}n^{\frac{3}{4}} \rfloor-1$ vertices, we are
guaranteed to have $m$ smallest $\omega_j$-neighbors. Suppose now
$x_1=\omega_j\cdot \rho_{y_{l_1},y_{r_1}}$ is the first connected
neighbor. Then we ``update'' $U_j(x_1)=U_j\dot\cup \{\rho_{
y_{l_1},y_{r_1}}\}$, $D_j(x_1)=D_j\dot\cup \{y_{l_1}\}$ and connect
the next $\omega_j$-neighbor via reversals contained in $N\setminus
U_j(x_1)$. Repeating this procedure until all smallest $m$
$\omega_j$-neighbors are explored, we obtain the set all connected
$\omega_j$-neighbors, $N[\omega_j]$. We then set
\begin{eqnarray*}
M_{j+1}&=& M_j\dot\cup N[\omega_j]\\
U_{j+1}&=&\cup_{x\in N[\omega_j]} U_j(x)\subset N\\
D_{j+1}&=&\cup_{x\in N[\omega_j]}D_j(x)
\subseteq\{\lfloor\frac{1}{2}n^{\frac{3}{4}}\rfloor+1,\ldots,n\}\\
L_{j+1}&=& L_j\setminus \{\omega_j\}\cup N[\omega_j].
\end{eqnarray*}
Note that each reversal is used at most once and reversals of the form
$\rho_{i,*}$ can only appear in one generation. \\
{\it Claim.} The above process generates a tree, that is each
$M_{j}$-element is
considered only once. \\
We prove the Claim by contradiction: assume the process generates a
cycle $\sigma_1\cdot\sigma_2\cdots\sigma_m\cdot\sigma_0=1$, where
$\sigma_i\in R_n$. Without loss of generality, we shall assume that
for reversals $\rho_{i,s}$, we always have $i\le s$ and consider
$j=\min\{h\mid \sigma_i=\rho_{h,s},0\le i\le m\}$. There are at most
two reversals among $\sigma_0,\sigma_1,\cdots,\sigma_m$ of the form
$\rho_{j,*}$. In case of only one such reversal, it is clear that
such a cycle cannot exist. Therefore we can, without loss of
generality, assume $\sigma_0=\rho_{j,a}$ and $\sigma_1=\rho_{j,b}$
where $a\ne b\in \mathbb{N}^+$. By construction, position $j$ is
never touched by the reversals $\sigma_2,\dots,\sigma_m$, whence we
arrive at the contradiction
\begin{equation*}
-b = (id\cdot\prod_{i=1}^{m}\sigma_{i})_j=
(id\cdot\sigma_0^{-1})_j=(id\cdot\sigma_0)_j=-a
\end{equation*}
and the Claim follows. Since $0<\delta<\frac{1}{4}$, we have for
sufficiently large $n$
\begin{equation*}
\mathbb{E}(N[w_j])=\frac{1+\epsilon_n}{\binom{n+1}{2}}\cdot
\left[\binom{n-\lfloor\frac{1}{2}n^{\frac{3}{4}}\rfloor+1}{2}
-(n-\lfloor\frac{1}{2}n^{\frac{3}{4}}\rfloor)
\lfloor\frac{1}{4}n^{\frac{3}{4}}\rfloor\right]>1.
\end{equation*}
Since we always connect only among the smallest $m$ neighbors, we
obtain, via the survival probability of the branching processes
$\mathcal{P}_P(\lambda)$ and
$\mathcal{P}_0(\lambda_n)$, depending
on whether we have $\epsilon_n=\epsilon>0$ or $\epsilon_n=o(1)$,
the following lower bound on $\mathbb{P}\left(\vert M_j\vert=\lfloor
\frac{1}{4}n^{\frac{3}{4}}\rfloor\mid \text{\rm for some $j$}
\right)$:
\begin{eqnarray*}
\mathbb{P}\left(\vert
M_j\vert=\lfloor\frac{1}{4}n^{\frac{3}{4}}\rfloor\mid \text{\rm for
some $j$} \right)\ge \wp(\epsilon_n)
\end{eqnarray*}
and the lemma follows.
\end{proof}
By choosing $k$ sufficiently large, we next enlarge the trees constructed
via Lemma~\ref{L:minigene} to subcomponents of arbitrary polynomial size,
which we call $k$-cells.
\begin{lemma}\label{L:expand2}
Suppose $k$ is arbitrary but fixed and $\theta_n\ge O(n^{\delta})$. Then
each $\Gamma_n$-vertex is contained in a $k$-cell, i.e.~a
$\Gamma_n$-subcomponent of size $\ge O(n^{\frac{3}{4}+k\delta})$
with probability at least
\begin{equation*}
\delta_k(\epsilon_n)=\wp(\epsilon_n)\, (1-e^{-\beta_k\theta_n}),
\quad\text{\it where } \beta_k>0.
\end{equation*}
\end{lemma}
\begin{proof}
Without loss of generality, we start the construction at the
identity. We set
\begin{equation*}
A_{m}=\left\{\rho_{l^{(m)},r^{(m)}}\in R_n \mid (m-1)\mu_n+1\le
l^{(m)}\le r^{(m)}\le m\mu_n\right\}
\end{equation*}
and write $w_{i}^{(h)}=\rho_{l_i^{(h)},r_i^{(h)}}\in A_h$. We
consider the branching process of Lemma~\ref{L:minigene} at $id$ and
denote the potentially generated tree of size
$\lfloor\frac{1}{4}n^{\frac{3}{4}}\rfloor$ by $T^1$. We consider the
r.v.
\begin{eqnarray*}
J_1=\left|\{w_{i}^{(1)}\in A_1\mid \exists x\in T^1;
\{x,x \cdot w_{i}^{(1)}\}\in \Gamma_n \}\right|.
\end{eqnarray*}
According to Lemma~\ref{L:minigene}, any two distinct
$J_1$-reversals connect distinct vertices
\begin{equation*}
\forall\, x,y\in T^1;\forall\, w_{i}^{(1)}\neq w_{r}^{(1)}\in
A_1;\quad x\cdot w_{i}^{(1)}\neq  y\cdot w_{r}^{(1)},
\end{equation*}
and the expected number of $J_1$-elements is given by
\begin{eqnarray}\label{E:exp1Y1}
\quad \mathbb{E}[J_1]=\binom{\mu_n+1}{2}\cdot
\left(1-\left(1-\frac{1+\epsilon_n}{\binom{n+1}{2}}\right)^{\lfloor
\frac{1}{4}n^{\frac{3}{4}}\rfloor}\right) \sim
\frac{1}{2}\mu_n^2\left(1-\exp(-(1+\epsilon_n)
\frac{1}{2}n^{-\frac{5}{4}})\right).
\end{eqnarray}
Chernoff's large deviation inequality
eq.~(\ref{E:cher}) \cite{Chernoff:52} implies that there exists
some constant $c_1>0$ such that
\begin{eqnarray*}
\mathbb{P}\left(J_1<\frac{1}{2}\mathbb{E}[J_1]\right)\le
\exp\left(-c_1\cdot \mathbb{E}[J_1]\right).
\end{eqnarray*}
We proceed by selecting the smallest element, $x_{j}^{(1)}$, from
the set $\{x \cdot w_{j}^{(1)}  \mid x\in T^1, \;w_j^{(1)}\in J_1\}$
and start the branching process of Lemma~\ref{L:minigene} at
$x_{j}^{(1)}$. As a result, we derive the subcomponent
$C_2(x_{j}^{(1)})$ of size $\lfloor
\frac{1}{4}n^{\frac{3}{4}}\rfloor$ with probability at least
$\wp(\epsilon_n)$. According to Lemma~\ref{L:minigene}, this
generation exclusively involves labels $j$ where
$j>\lfloor\frac{1}{2}n^{\frac{3}{4}}\rfloor$.
Therefore, since any two smallest elements $x_{j_1}^{(1)}$ and
$x_{j_2}^{(1)}$ differ in at least one coordinate, $h$, for $1\le
h\le \mu_n$, which is not touched by the branching process of
Lemma~\ref{L:minigene}, we have
\begin{equation*}
C_2(x_{j_1}^{(1)})\cap C_2(x_{j_2}^{(1)})=\varnothing.
\end{equation*}
Let $X_1$ be the r.v.~counting the number of these new
$\Gamma_n$-subcomponents. In view of eq.~(\ref{E:exp1Y1}), we obtain
\begin{equation*}
\mathbb{E}[X_1] \ge  \wp(\epsilon_n)\cdot \mathbb{E}[J_1]\sim
\wp(\epsilon_n)\cdot\frac{1}{2}\mu_n^2
\left(1-\exp(-(1+\epsilon_n)\frac{1}{2}n^{-\frac{5}{4}})\right)\triangleq
\theta_n \ge O(n^{\delta}).
\end{equation*}
Again, using the large deviation inequality, eq.~(\ref{E:cher}), we conclude
that there exists some $\beta_1>0$
such that
\begin{equation*}
\mathbb{P}(X_1<\frac{1}{2}\theta_n)\le
\exp(-\beta_1\theta_n).
\end{equation*}
The union of all the
$C_2(x_{j}^{(1)})$-subcomponents with $T^1$ forms a
$\Gamma(B_n,R_n)$-subcomponent, $T^2$, and we have
\begin{equation*}
\mathbb{P}\left(\vert T^2\vert < \lfloor
\frac{1}{4}n^{\frac{3}{4}}\rfloor \cdot \frac{1}{2}\theta_n\right)
\le \exp(-\beta_1\theta_n).
\end{equation*}
We proceed by induction:\\
{\it Claim:} For each $1\le i\le k$, there exists some constant
$\beta_{i-1}>0$ and a
$\Gamma(B_n,R_n)$-subcomponent $T^i$ such that
\begin{eqnarray*}
\mathbb{P}(\vert T^i\vert <\lfloor \frac{1}{4}n^{\frac{3}{4}}\rfloor
\cdot
\frac{1}{2^{i-1}}\theta_n^{i-1})\le
\exp(-\beta_{i-1}\theta_n).
\end{eqnarray*}
We have already established the induction basis. As for the
induction step, let us assume the Claim holds for $i<k$ and let
$C_i(\alpha)$ denote a subcomponent generated by the branching
process of Lemma~\ref{L:minigene} in the $i$-th step at $\alpha$.
We consider
the $w_{r}^{(i+1)}\ne w_{a}^{(i+1)}\in A_{i+1}$ and
\begin{eqnarray*}
J_{i+1}=\{w_{r}^{(i+1)}\in A_{i+1}\mid \exists x\in C_i(\alpha);\;
\{x,x \cdot w_{r}^{(i+1)}\}\in \Gamma_n
\}.
\end{eqnarray*}
At the minimal elements, $x^{\alpha}_{r}$ of $\{x \cdot
w_{r}^{(i+1)} \mid x\in C_i(\alpha),\, w_{r}^{(i+1)}\in J_{i+1}\}$,
we initiate the branching process of Lemma~\ref{L:minigene}. The
process generates subcomponents $C_{i+1}(x^\alpha_{r})$ of size
$\lfloor \frac{1}{4}n^{\frac{3}{4}}\rfloor$ with probability $\ge
\wp(\epsilon_n)$. By construction, any two of these are mutually
disjoint and let $X_{i+1}$ be the r.v.~counting their number. We
derive setting $q_n=\lfloor \frac{1}{4}n^{\frac{3}{4}}\rfloor$
\begin{eqnarray*}
\mathbb{P}\left(\vert T^{i+1}\vert < q_n
\frac{1}{2^{i}}\theta_n^{i}\right)
& \le & \underbrace{\mathbb{P}\left( \vert T^i\vert <q_n
\frac{1}{2^{i-1}}\theta_n^{i-1}\right)
}_{\text{\rm failure at step $i$}} \ +
\underbrace{ \mathbb{P}\left(\vert
T^{i+1}\vert < q_n\frac{1}{2^{i}}\theta_n^{i}\mbox{ and }\vert
T^i\vert \ge q_n\frac{1}{2^{i-1}}\theta_n^{i-1} \right)
}_{\text{\rm failure at step $i+1$ conditional to $\vert T^i\vert
\ge
q_n\frac{1}{2^{i-1}}\theta_n^{i-1}$}}\\
& \le &
\underbrace{e^{-\beta_{i-1}\,\theta_n}}_{\text{\rm
induction hypothesis}} +
\underbrace{e^{-\beta
\,\theta_n^{i}}}_{\text{\rm large deviation results}}\cdot
(1-e^{-\beta_{i-1}\, \theta_n})\, , \quad \beta_{i-1}>0 \\
& \le & e^{-\beta_{i}\,\theta_n}
\end{eqnarray*}
and the Claim follows.\\
Therefore each $\Gamma_n$-vertex is contained in a subcomponent of
size at least $O(n^{\frac{3}{4}+k\delta})$, with probability at
least $\wp(\epsilon_n)(1-\exp(-\beta_k\theta_n))$ and the lemma is
proved.
\end{proof}

We will call a subcomponent constructed in Lemma~\ref{L:expand2} a
$k$-cell or simply a cell.

\section{Small components}\label{S:4}
Let $\Gamma_{n,k}$ denote the set of $\Gamma_n$-vertices contained
in components of size $\ge O(n^{\frac{3}{4}+k\delta})$ for some
$0<\delta<\frac{1}{4}$. In this section we prove that $\vert
\Gamma_{n,k}\vert$ is a.s.~$ \sim \wp(\epsilon_n)\cdot 2^n\cdot n!$.
In analogy to Lemma~$3$ of \cite{Reidys:08rand} we first observe
that the number of vertices, contained in $\Gamma_n$-components of
size $<c_k\,n^{\frac{3}{4}+k\delta}$, is sharply concentrated. The
concentration reduces the problem to a computation of expectation
values. It follows from considering the indicator r.vs. of pairs
$(C,v)$ where $C$ is a component and $v\in C$ and to estimate their
correlation. Since the components in question are small, no
``critical'' correlation terms arise.
\begin{lemma}\label{L:comple1}\cite{Reidys:08rand}
Let $\omega_n=\vert \Gamma_n\setminus \Gamma_{n,k}\vert$ and
$\lambda_n=\frac{1+\epsilon_n}{\binom{n+1}{2}}$, where
$n^{-\frac{1}{4}+\delta}\le \epsilon_n \le \lambda$, for some
$\lambda>0$. Then we have
\begin{equation*}
\mathbb{P}\left(\mid \omega_n -\mathbb{E}[\omega_n]\mid \, \ge \,
\frac{1}{n} \mathbb{E}[\omega_n]\right)=o(1) .
\end{equation*}
\end{lemma}
With the help of Lemma~\ref{L:comple1}, we are in position to compute the
size of $\Gamma_{n,k}$.
\begin{lemma}\label{L:size1}
Let $\lambda_n=\frac{1+\epsilon_n}{\binom{n+1}{2}}$, where
$n^{-\frac{1}{4}+\delta}\le \epsilon_n\le \epsilon <1$ and suppose
$k\in\mathbb{N}$ is sufficiently large. Then
\begin{equation*}
\vert \Gamma_{n,k}\vert \sim \wp(\epsilon_n)\cdot 2^n\cdot n! \
\qquad \text{\it a.s.~.}
\end{equation*}
\end{lemma}
\begin{proof}
First we prove for any $n^{-\frac{1}{4}+\delta}\le \epsilon_n\le
\lambda$, where $\lambda>0$
\begin{equation*}
(1-o(1))\wp(\epsilon_n)\cdot 2^n\cdot n!\le \vert \Gamma_{n,k}\vert
\qquad \text{\rm  a.s.}
\end{equation*}
By assumption we have
$$
\mathbb{E}[\omega_n]\le (1-\delta_{k}(\epsilon_n))\cdot 2^n\cdot n!.
$$
In view of Lemma~\ref{L:comple1}, we derive
\begin{equation*}
\omega_n < \left(1+O(\frac{1}{n})\right) \,\mathbb{E}[\omega_n]
<\left(1-\delta_k(\epsilon_n)+O(\frac{1}{n})\right)\cdot 2^n\cdot n!
\quad \text{\rm a.s.,}
\end{equation*}
whence
\begin{eqnarray*}
\vert\Gamma_{n,k}\vert
\ge\left(\delta_k(\epsilon_n)-O(\frac{1}{n})\right)\vert
\Gamma_n\vert=(1-o(1))\wp(\epsilon_n)\cdot 2^n\cdot n!\quad
\text{\rm a.s..}
\end{eqnarray*}
Next we prove for $n^{-\frac{1}{4}+\delta}\le \epsilon_n < 1$
\begin{equation*}
\vert\Gamma_{n,k}\vert\le (1+o(1))\wp(\epsilon_n)\cdot 2^n\cdot
n!\qquad \text{\rm a.s.}
\end{equation*}
For this purpose we consider the branching process on a
$\binom{n+1}{2}$-regular rooted
tree $T_{r^*}$ where the r.v. $\xi_r^*$ of the rooted vertex $r^*$
is $B(\binom{n+1}{2},\lambda_n)$ distributed while the r.v. of any
other vertex $r$ has the distribution
$B(\binom{n+1}{2}-1,\lambda_n)$. Let $C_{r^*}$ denote the component
generated by such a branching process. Bollob\'{a}s {\it et al.}
\cite{Bollobas:92} showed that
\begin{eqnarray}\label{E:ts}
\quad
\mathbb{P}(\vert C_{r^*}\vert=i)=(1+o(1))\cdot\frac{(\lambda_n\cdot
(\binom{n+1}{2}-1))^{i-1}}{i\sqrt{2\pi
i}}\left[\frac{(\binom{n+1}{2}-1)(1-\lambda_n)}{\binom{n+1}{2}-2}
\right]^{(\binom{n+1}{2}-2)i+2},
\end{eqnarray}
where $i=i(n)\rightarrow\infty$ as $n\rightarrow\infty$. The key
observation is an inequality \cite{Bollobas:92}, relating this
process with the construction of a spanning component of a
$\Gamma_n$-component at vertex $r$,
\begin{equation}\label{E:retrees}
\mathbb{P}\left(\vert C_{r^*}\vert \le m\right)\le
\mathbb{P}\left(\vert C_{r}\vert \le m\right).
\end{equation}
Eq.~(\ref{E:retrees}) follows immediately from the observation that
during the generation of a spanning component, there are for each
vertex at most $(\binom{n+1}{2}-1)$ neighbors that are not in the
component, while in $T_{r^*}$ there exist exactly
$(\binom{n+1}{2}-1)$ new neighbors. Suppose now $k$ is sufficiently
large, satisfying $k\delta+\frac{3}{4}\ge 3\frac{1}{3}$. Then
$n^2\le c_k\cdot n^{k\delta+\frac{3}{4}}$ for sufficiently large
$n$, i.e.
\begin{equation*}
\mathbb{P}(c_k\cdot n^{k\delta+\frac{3}{4}}\le \vert
C_{r^*}\vert<\infty) \le \mathbb{P}(n^2\le \vert
C_{r^*}\vert<\infty).
\end{equation*}
This probability can be estimated as follows
\begin{eqnarray*}
\mathbb{P}(n^2\le \vert C_{r^*}\vert < \infty) &=& \sum_{i\ge
n^2}\mathbb{P}(\vert
C_{r^*}\vert=i)\\
& \le &  \sum_{i\ge n^2} (1 +
o(1))\cdot \frac{(\lambda_n\cdot
(\binom{n+1}{2}-1))^{i-1}}{i\sqrt{2\pi
i}}\left[\frac{(\binom{n+1}{2}-1)(1-\lambda_n)}
{(\binom{n+1}{2}-2)}\right]^{(\binom{n+1}{2}-2)i+2}\\
& \le & \sum_{i\ge n^2}\left[(1+\epsilon_n)e^{-\epsilon_n}
\right]^i\le \sum_{i\ge n^2} c(\epsilon)^i =o(e^{-n\ln(2n)}),
\end{eqnarray*}
where $0<c(\epsilon)<1$. We accordingly derive
\begin{eqnarray}
\mathbb{P}(\vert C_{r^*}\vert<c_k\cdot n^{k\delta+\frac{3}{4}}) &= &
\mathbb{P}(\vert C_{r^*}\vert <\infty)-
\mathbb{P}(c_k\cdot n^{k\delta+\frac{3}{4}} \le \vert C_{r^*}\vert
< \infty)\\
&\ge& \label{E:est}
(1-\underbrace{(1+o(1))\wp(\epsilon_n)}_{=\pi_0})-o(e^{-n\ln(2n)}),
\end{eqnarray}
where $\pi_0$ denotes the survival probability of the branching
process on $T_{r^*}$, see Corollary~\ref{C:1}. From
eq.~(\ref{E:retrees}) and eq.~(\ref{E:est}) we immediately obtain,
taking the expectation
\begin{equation*}
\mathbb{E}[\omega_n]\ge (1-\wp(\epsilon_n))\cdot 2^n\cdot n!+o(1).
\end{equation*}
Lemma~\ref{L:comple1} accordingly implies
\begin{equation*}\label{E:jk1}
(1-\wp(\epsilon_n)-O(\frac{1}{n}))\cdot 2^n\cdot n! \le \omega_n
\qquad \text{\rm a.s.,}
\end{equation*}
whence the lemma.
\end{proof}

\begin{lemma}\label{L:1-e}
Let $\lambda_n=\frac{1-\epsilon}{\binom{n+1}{2}}$, where
$0<\epsilon<1$, then a.s.~$\Gamma_n$ contains no component
larger than $O(n\ln(n))$.
\end{lemma}
\begin{proof}
We show that there exists some $\kappa>0$ such that $\vert
C_n^{(1)}\vert\le\kappa\cdot n\ln(n)$ a.s.. For this purpose we
study the probability that each vertex $r$ is contained in a
component of size $>\kappa\cdot n\ln(n)$. As in Lemma~\ref{L:size1},
let $C_r$ be the component containing vertex $r$ in $\Gamma_n$ and
let $C_{r*}$ denote the component in the
$\binom{n+1}{2}$-regular tree
$T_{r*}$, rooted in $r^*$. The key observation is
eq.~(\ref{E:retrees}),
\begin{equation*}
\mathbb{P}\left(\vert C_{r^*}\vert \le m\right)\le
\mathbb{P}\left(\vert C_{r}\vert \le m\right),
\end{equation*}
which implies
\begin{eqnarray*}
\mathbb{P}\left(\vert C_r \vert>\kappa\cdot n\ln(n)\right) &\le
&1-\mathbb{P}\left(\vert C_{r^*}
\vert\le\kappa\cdot n\ln(n)\right)\\
&=&\sum_{i>\kappa\cdot n\ln(n)}\mathbb{P}\left(\vert C_{r^*}
\vert=i\right).
\end{eqnarray*}
Let $X_\kappa$ denote the r.v counting the number of vertices
in components of size $>\kappa\cdot n\ln(n)$. In view of eq.~(\ref{E:ts})
and $\lambda_n=\frac{1-\epsilon}{\binom{n+1}{2}}$, we derive
\begin{eqnarray*}
\mathbb{E}(X_{\kappa})\le 2^n\cdot
n!\cdot\mathbb{P}\left(\vert C_r \vert>\kappa\cdot n\ln(n)\right) &\le &
\frac{2^n\cdot n!}{1-\epsilon}\cdot\sum_{i>\kappa\cdot
n\ln(n)}i^{-\frac{3}{2}}((1-\epsilon)\cdot e^{\epsilon})^i.
\end{eqnarray*}
Since $\Lambda_\epsilon=(1-\epsilon)\cdot e^{\epsilon}<1$, for any
$\epsilon >0$, we arrive at
\begin{eqnarray*}
\mathbb{E}(X_{\kappa})
&\le &\frac{2^n\cdot n!}{1-\epsilon}\cdot (\kappa\cdot
n\ln(n))^{-\frac{3}{2}}\Lambda_\epsilon^{\kappa\cdot
n\ln(n)}\cdot \frac{1}{1-\Lambda_\epsilon}=o(1).
\end{eqnarray*}
Thus, choosing $\kappa$ sufficiently large, we can conclude that
a.s. $\Gamma_n$ contains no components of size $>\kappa\cdot n\ln(n)$.
\end{proof}

\section{Density and splits}\label{S:dense}

\begin{lemma}\label{L:dense-01}
$\Gamma_{n,k}$ is a.s.~dense in $B_n$.
\end{lemma}
\begin{proof}
We consider
\begin{equation*}
A_{k+1}=\left\{\rho_{l_j^{(k+1)},r_j^{(k+1)}}\in R_n \mid
k\mu_n+1\le l_j^{(k+1)}\le
r_j^{(k+1)}\le\lfloor\frac{1}{2}n^{\frac{3}{4}}\rfloor\right\}
\end{equation*}
Let $w_{j}^{(k+1)}=\rho_{l_j^{(k+1)},r_j^{(k+1)}}$ and recall that
$\ell_n=\lfloor \frac{k}{2(k+1)} n^{\frac{3}{4}}\rfloor$. We set
\begin{equation*}
d^{(k+1)}(v)=\{v\cdot w_{i}^{(k+1)}\mid 1\le i\le
\binom{\ell_n+1}{2}\}.
\end{equation*}
Clearly,
\begin{eqnarray*}
\vert d^{(k+1)}(v)\vert &=& \binom{\ell_n+1}{2}
\sim \frac{1}{2}\left[\frac{k}{2(k+1)}\right]^2 \cdot
n^{\frac{3}{2}}\cdot(1+o(1)).
\end{eqnarray*}
Let $\Delta_k=\left[\frac{k}{2(k+1)}\right]^{2}/2$ and $Z(v)$ be the
r.v.\ counting the number of vertices contained in the set
$d^{(k+1)}(v)\cap\Gamma_{n,k}$, whose subcomponents are constructed
in Lemma~\ref{L:expand2}. We immediately observe
\begin{equation*}
\mathbb{E}(Z(v))\ge \delta_k(\epsilon_n)\cdot\vert d^{(k+1)}(v) \vert
\sim
\Delta_k\,n^{\frac{3}{2}}\cdot\wp(\epsilon_n)(1-e^{-\beta_k\theta_n})
\ge \Delta_k\cdot n^{\frac{5}{4}+\delta}.
\end{equation*}
Since the construction of the Lemma~\ref{L:expand2}-subcomponents did not
involve any elements contained in
$[k\mu_n+1,\lfloor\frac{1}{2}n^{\frac{3}{4}}\rfloor]$, any two such
subcomponents are vertex-disjoint. Therefore the r.v.\ $Z(v)$ is a
sum of {\it independent} indicator r.vs.~and Chernoff's large
deviation inequality \cite{Chernoff:52} implies
\begin{equation*}
\mathbb{P} \left( Z(v) <
\frac{1}{2}\,\Delta_k\cdot n^{\frac{5}{4}+\delta}\right)\le
\exp(-\kappa n^{\frac{5}{4}+\delta}) \quad \text{\rm for some
$\kappa>0$.}
\end{equation*}
We conclude from this that the expected number of $B_n$-vertices with the
property $Z(v) < \frac{1}{2}\,\Delta_k\,  n^{\frac{5}{4}+\delta}$ is
tending to zero, whence the lemma.
\end{proof}
Next we show that there exist many vertex disjoint paths between
$\Gamma_{n,k}$-splits of sufficiently large size. The proof is
analogous to Lemma~$7$ in \cite{Reidys:08rand}. We remark that
Lemma~\ref{L:split1} does not use an isoperimetric inequality
\cite{Harper:66b}. It only employs a generic estimate of the vertex
boundary in Cayley graphs due to Aldous \cite{Aldous:87,Babai:91b}.
\begin{lemma}\label{L:split1}
Let $(S,T)$ be a vertex-split of $\Gamma_{n,k}$ with the properties
\begin{equation}\label{E:prop}
\exists\,0< \rho_0\le \rho_1<1;\quad 2^n\cdot(n-2)!\le \vert S\vert
= \rho_0 \vert \Gamma_{n,k}\vert  \quad \text{\rm and}\quad
2^n\cdot(n-2)!\le \vert T \vert = \rho_1 \vert \Gamma_{n,k}\vert .
\end{equation}
Then there exists some $c>0$ such that a.s.~$S$ is connected to $T$
in $\Gamma(B_n,R_n)$ via at least
\begin{equation*}
c\cdot\frac{2^n\cdot(n-3)!}{\binom{n+1}{2}^3}
\end{equation*}
edge disjoint (independent) paths of length $\le 3$.
\end{lemma}
\begin{proof}
We distinguish the cases $\vert \text{\sf B}(S,1)\vert\le
\frac{2}{3}\cdot 2^n\cdot n!$ and $\vert\text{\sf B}(S,1)\vert>
\frac{2}{3}\cdot 2^n\cdot n!$. In the former case, we employ the
generic estimate of vertex boundaries in Cayley graphs \cite{Aldous:87},
\begin{equation}\label{E:gen}
\vert {\sf d}(A) \vert \ge \frac{1}{{\rm
diam}(\Gamma(B_n,R_n))}\cdot\vert A \vert \left(1-\frac{\vert
A\vert}{2^n\cdot n!}\right).
\end{equation}
In view of eq.~(\ref{E:prop}) and ${\rm diam}(\Gamma(B_n,R_n))=n+1$
\cite{Konstantinova}, eq.~(\ref{E:gen}) implies
\begin{equation*}
\exists \,d_1>0;\quad \vert {\sf d}({\sf B}(S,1))\vert
\ge\frac{d_1}{n+1}\cdot\vert{\sf B}(S,1)\vert \ge d_1\cdot 2^{n-1}
\cdot (n-3)! .
\end{equation*}
According to Lemma~\ref{L:dense-01}, a.s.~all signed permutations are
within distance $1$ to some $\Gamma_{n,k}$-vertex, whence
\begin{equation*}
\vert \text{\sf d}(\text{\sf B}(S,1))\cap \text{\sf B}(T,1)\vert\ge
d_1\cdot 2^{n-1}\cdot (n-3)! \quad \text{\rm a.s..}
\end{equation*}
Let $\beta_1\in \text{\sf d}(\text{\sf B}(S,1))\cap \text{\sf
B}(T,1)$ and set
\begin{equation*}
T^*=\{\alpha_1 \in \text{\sf d}(S)\mid
d(\alpha_1,\beta_1)=1,\text{\rm for some $\beta_1\in \text{\sf
B}(T,1)$}\}.
\end{equation*}
Evidently, at most $\binom{n+1}{2}$ elements in $\text{\sf d}(S)$
can be connected to the same $\beta_1$, whence there are a.s.~at
least
$$
 d_{1}\cdot\frac{2^{n-1}\cdot
(n-3)!}{\binom{n+1}{2}}
$$
edge disjoint paths connecting ${\sf d}(S)$ to ${\sf B}(T,1)$. Let
furthermore $T_1\subset T^*$ be some maximal set such that any pair
of $T_1$-vertices $(\beta_1,\beta_1')$ has distance
$d(\beta_1,\beta_1')> 2$. Then $\vert T_1\vert > \vert
T^*\vert/2(\binom{n+1}{2}+1)^2$
since $\vert \text{\sf B}(v,2)\vert<
{2(\binom{n+1}{2}+1)^2}$. By construction, any two of the paths from
$S$ to $T_1\subset \text{\sf d}(S)$ are edge disjoint and
accordingly there are a.s.~at least
\begin{equation*}
d_{1}\cdot \frac{ 2^{n-1}\cdot(n-3)!
}{{
2}\binom{n+1}{2}(\binom{n+1}{2}+1)^2} \sim
c\cdot\frac{2^n\cdot(n-3)!}{\binom{n+1}{2}^3}, \quad \text{\rm
where}\ c>0
\end{equation*}
edge disjoint paths of length $2$ or $3$ connecting $S$ and
$T$. \\
It remains to study the case $\vert \text{\sf B}(S,1)\vert >
\frac{2}{3}\cdot 2^n\cdot n!$. By construction both: $S$ and $T$
satisfy eq.~(\ref{E:prop}), whence we can, without loss of
generality assume that also $\vert \text{\sf B}(T,1)\vert >
\frac{2}{3}\cdot 2^n\cdot n!$ holds. But then
$$
\vert \text{\sf B}(S,1)\cap \text{\sf B}(T,1)\vert >
\frac{1}{3}\cdot 2^n\cdot n!,
$$
we have a.s~at least $ \frac{1}{3}\cdot 2^n\cdot n!$ edge disjoint
paths of length $\le 2$ connecting $S$ and $T$.
\end{proof}

